\documentclass[12pt,oneside,british,a4wide]{amsart}
\usepackage[T1]{fontenc}
\usepackage[latin9]{inputenc}
\usepackage{amstext}
\usepackage{amsthm}
\usepackage{amssymb}
\usepackage{enumitem}

\makeatletter
\numberwithin{equation}{section}
\numberwithin{figure}{section}

\usepackage{amsmath, wrapfig}
\usepackage{dsfont, a4wide, amsthm, amssymb, amsfonts, graphicx}
\usepackage{fancyhdr, xspace, psfrag, setspace, supertabular, color, cite}
\usepackage[pagebackref=true]{hyperref}

\newtheorem{theorem}{Theorem}[section]
\newtheorem{proposition}[theorem]{Proposition}
\newtheorem{lemma}[theorem]{Lemma}
\newtheorem{claim}[theorem]{Claim}
\newtheorem{corollary}[theorem]{Corollary}

\newtheorem{conjecture}[theorem]{Conjecture}

\newtheorem{fact}[theorem]{Fact}

\theoremstyle{definition}
\newtheorem{definition}[theorem]{Definition}

\def\showauthornotes{1}

\ifnum\showauthornotes=1
\newcommand{\Authornote}[2]{{\sf\small\color{red}{[#1: #2]}}}
\else
\newcommand{\Authornote}[2]{}
\fi

\def\E{\mathbb{E}}

\def\C{\mathbb{C}}

\def\F{\mathbb{F}}

\def\VC{\mathrm{VC}}

\def\Err{\mathrm{Err}}

\def\e{\epsilon}

\newcommand{\abar}{\bar{a}}
\newcommand{\bbar}{\bar{b}}

\newcommand{\calL}{\mathcal{L}}
\newcommand{\calF}{\mathcal{F}}

\newcommand{\calC}{\mathcal{C}}

\newcommand{\calB}{\mathcal{B}}
\newcommand{\calD}{\mathcal{D}}
\newcommand{\calP}{\mathcal{P}}

\newcommand{\calG}{\mathcal{G}}

\newcommand{\calQ}{\mathcal{Q}}

\newcommand{\calS}{\mathcal{S}}

\newcommand{\codim}{{\mathrm{codim}}}

\newcommand{\At}{{\mathrm{At}}}

\newcommand{\Bool}{{\mathrm{Bool}}}
\newcommand{\Sym}{{\mathrm{Sym}}}
\newcommand{\ind}{{\mathrm{ind}}}

\newenvironment{proofof}[1]{\indent{\scshape Proof of #1}:~~}{\qed}

\makeatother

\usepackage{babel}
\begin{document}

\title{A localized packing lemma for sets of bounded $\VC_2$-dimension}

\author{C. Terry}

\author{J. Wolf}

\address{Department of Mathematics, Statistics, and Computer Science, University of Illinois Chicago, Chicago IL 60607, USA}

\email{caterry@uic.edu}

\address{Department of Pure Mathematics and Mathematical Statistics, Centre for Mathematical Sciences, Wilberforce Road, Cambridge CB3 0WB, UK}

\email{julia.wolf@dpmms.cam.ac.uk}

\date{}

\begin{abstract} In the setting of abelian groups, Haussler's packing lemma implies that any subset of bounded $\VC$-dimension can be approximated by one of a small number of representative translates. While this exact statement is provably false under the weaker hypothesis that $A$ has bounded $\VC_2$-dimension, we show that in $\F_p^n$ a similar statement does hold after localization to cosets of subspaces depending on the individual translates. 

More precisely, given $\e>0$ and $A\subseteq \F_p^n$ with $\VC_2(A)<k$, we show that there exist $r\leq \e^{-O_{k,p}(1)}$ and $g_1,\ldots, g_r\in \F_p^n$ with the property that for every $x\in \F_p^n$, there is a partition $\calC_x$ of $\F_p^n$ into cosets of (possibly distinct) subgroups, each of index at most $\exp(\e^{-O_{k,p}(1)})$, and a subcollection $\Omega_x\subseteq \calC_x$ whose union has size at least $(1-\e)|\F_p^n|$, such that on each $C\in \Omega_x$, the translate $A-x$ looks like one of $A-g_1,\ldots, A-g_r$, in the sense that $|((A-x)\Delta (A-g_{i_{x,C}}))\cap C|<\e|C|$  for some $1\leq i_{x,C}\leq r$.  We also prove a second version of this result in which each partition $\calC_x$ is the set of cosets of a single subgroup $H_x$ of quasipolynomial index but where the bound on $r$ is a quadruple exponential.
\end{abstract}

\maketitle


\section{Introduction}

Haussler's packing lemma \cite{Haussler.1995} is a foundational result in $\VC$-theory which has had countless applications across combinatorics, statistics and computer science.  In the setting of groups, Alon, Fox and Zhao in \cite{Alon.2019k3} observed that the packing lemma has the following consequence for a set whose translates form a set system of bounded $\VC$-dimension (see Definition \ref{def:vc}). 
 
\begin{corollary}[Immediate consequence of Haussler packing]\label{thm:AFZ}
Let $k\geq 1$, let $\delta \in (0,1]$, let $G$ be a finite abelian group, and let $A\subseteq G$ be a subset of $\VC$-dimension at most $k$. Then there are  $t\leq (30/\delta)^k$ and $g_1,\ldots, g_t\in G$ with the property that for all $x\in G$, there is $i\in [t]$ such that $|(A-x)\Delta (A-g_i)|\leq \delta|G|$.
\end{corollary}

Roughly speaking, Corollary \ref{thm:AFZ}  tells us that every translate of $A$ looks approximately like one of a very small number of ``representative'' translates. We first show that this statement can fail spectacularly under the weaker hypothesis of bounded $\VC_2$-dimension (see Section \ref{ss:vc}). For example, the set 
\begin{align}\label{al:atom}
A=\{x\in \F_p^n: x\cdot x=0\}
\end{align}
has $\VC_2$-dimension $1$ (see \cite[Lemma 3.10]{Terry.2025}, which originally appeared in an early version of \cite{Terry.2021}). However, we shall see in Proposition \ref{prop:atom} that, assuming $p$ is odd and $\delta$ is sufficiently small, any collection $g_1,\ldots, g_t$ satisfying the conclusions of Corollary \ref{thm:AFZ} for this set $A$ must have $t\geq n$. Thus no literal analogue of Corollary \ref{thm:AFZ} can be proved with a bound depending only on $\VC_2(A)$ and $\delta$.   On the other hand, we will show that the set in (\ref{al:atom}) does satisfy a packing statement \emph{after localizing to cosets of subspaces depending on individual translates}.   
 
\begin{proposition}[Localized packing for a quadratic atom]\label{exa:quadat}
Let $p>2$ be a prime and let $n\geq 3$ be an integer. Let $G=\F_p^n$ and let $A=\{x\in G: x\cdot x=0\}$.  Then there exist elements $g_1,\dots, g_{2p} \in G$ and for every $x\in G$, a subspace $H_x\leqslant G$ of codimension at most 3, such that for every coset $C$ of $H_x$, there is some $i_{x,C}\in [2p]$ satisfying
\[(A-x)\cap C=(A-g_{i_{x,C}})\cap C.\]
\end{proposition}

 The conclusions of Proposition \ref{exa:quadat} are substantially weaker than those in Corollary \ref{thm:AFZ}. In particular, for each translate $A-x$, Proposition \ref{exa:quadat} introduces a partition of $G$ into cosets of a subgroup $H_x$, depending on $x$, so that $A-x$ looks like some  element from $\{A-g_1,\ldots, A-g_{2p}\}$ on each coset of $H_x$.  On different cosets of $H_x$, $A-x$ might look like different elements of $\{A-g_1,\ldots, A-g_{2p}\}$, but there is a single, global list of ``representative'' translates. For this reason we consider Proposition \ref{exa:quadat} a ``localized packing lemma''.
 
The goal of this paper is to prove that all subsets of $\F_p^n$ of bounded $\VC_2$-dimension satisfy a similar localized packing lemma, with efficient bounds on both the number of translates required and the codimension of the subgroups involved.  Our first such result is stated below.

\begin{theorem}[Localized packing with polynomially many translates]\label{thm:packingintro}
For all primes $p>2$, integers $k\geq 1$, and $0<\e<1/2$, there exist integers 
$$
 R\leq  \e^{-O_{k,p}(1)}\text{ and }D\leq \e^{-O_{k,p}(1)}
$$
with the following property.  

Let $G=\F_p^n$ and suppose that $A\subseteq G$ satisfies $\VC_2(A)< k$.  Then there exist $r\leq R$ and elements $g_1,\ldots,g_r\in G$ such that for every $x\in G$, there is a partition $\mathcal C_x$ of $G$ and a subcollection $\Omega_x\subseteq\mathcal C_x$ such that  
\begin{enumerate}[label=\normalfont(\roman*)] 
\item  the elements of $\calC_x$ are cosets of (possibly distinct) subgroups of $G$, each of codimension at most $D$;
\item $|\bigcup_{C\in\Omega_x}C|\geq (1-\epsilon)|G|$;
\item for every $C\in\Omega_x$, there is $i_{x,C}\in[r]$ such that $|((A-x)\triangle(A-g_{i_{x,C}})\bigr)\cap C |<\epsilon |C|$. 
\end{enumerate}
\end{theorem}

 We prove Theorem \ref{thm:packingintro} using a variation on a sufficient condition for large $\VC_2$-dimension developed in  \cite{Terry.2021}, and applied in \cite{Terry.2025quo}, along with a variation on the $U^2$ cylinder regularity lemma from \cite{Fox.2022}. 
 
A more sophisticated version of Theorem \ref{thm:packingintro} will be used in forthcoming work \cite{Terry.2027} of the authors, combined with Bogolyubov-type arguments, to give a new proof of the structure theorem for sets of bounded $\VC_2$-dimension first proved by the authors in \cite{Terry.2021} (now \cite{Terry.2025}). However, we believe that the packing result is of independent interest. 
 
It is reasonable to ask whether the bound on the codimension of the cosets appearing in Theorem \ref{thm:packingintro} is best possible. In particular, in \cite{Terry.2025quo} the authors proved that the quadratic partition size required for a full structure theorem for sets of bounded $\VC_2$-dimension was $O(\epsilon^{-k})$, which would necessitate a logarithmic bound on $D$.

Towards this end, we also prove a version of Theorem \ref{thm:packingintro} in which $D$ is sub-polynomial, at the cost of worsening the bound on the number of translates required. This second version also upgrades each $\calC_x$ to a partition consisting of cosets of a single subgroup $H_x$ (although the subgroup $H_x$ still varies with $x$). 

\begin{theorem}[Localized packing with quasipolynomial index]\label{thm:packingintro2}
For all primes $p>2$, integers $k\geq 1$, and $0<\e<1$, there exist integers 
$$
 R\leq  \exp_p(\exp_p(\exp_p(\exp_p(\e^{-O_{k,p}(1)}))))\text{ and }D\leq 8 \log^2_p((\e/4)^{-2k-o(1)})
$$
with the following property. 

Let $G=\F_p^n$ with $n$ sufficiently large, and suppose that $A\subseteq G$ satisfies $\VC_2(A)< k$.  Then there exist $r\leq R$ and elements $g_1,\ldots,g_r\in G$ such that for every $x\in G$, there is a partition $\mathcal C_x$ of $G$ and a subcollection $\Omega_x\subseteq\mathcal C_x$ such that  
\begin{enumerate}[label=\normalfont(\roman*)] 
\item  $\calC_x=G/H_x$ for some subgroup $H_x$ of codimension at most $D$;
\item $|\bigcup_{C\in\Omega_x}C|\geq (1-\epsilon)|G|$;
\item for every $C\in\Omega_x$, there is $j_{x,C}\in[r]$ such that $|((A-x)\triangle(A-g_{j_{x,C}}))\cap C |<\epsilon |C|$. 
\end{enumerate}
Here $o(1)$ denotes a quantity that tends to 0 as $\e$ tends to zero, at a rate that may depend on $p$ and $k$.
\end{theorem}
 
 We deduce Theorem \ref{thm:packingintro2} from a full quadratic structure theorem for sets of bounded $\VC_2$-dimension, in particular, the fact that any set of bounded $\VC_2$-dimension is approximately equal to a union of atoms of a quadratic factor of bounded complexity.  Such a result was first proved in \cite{Terry.2021} (now \cite{Terry.2025}), and was subsequently refined quantitatively,  first in \cite{Terry.2025quo}, and again in \cite{Sheats.2025}. Our proof of Theorem \ref{thm:packingintro2} shows that any such structure theorem gives rise to a localized packing lemma, with the specific bounds in Theorem \ref{thm:packingintro2} resulting from \cite[Theorem 1.3]{Sheats.2025}.  

Interestingly, we know of no way to simultaneously obtain the polynomial bound on $R$ in Theorem \ref{thm:packingintro}, and the sub-polynomial bound on $D$ in Theorem \ref{thm:packingintro2}.  We conjecture that such a version is achievable.

\begin{conjecture}
One can obtain the conclusion of Theorem \ref{thm:packingintro2} with $R\leq \e^{-O_{k,p}(1)}$ and $D\leq \log^2_p(O_{k,p}(\e^{-1}))$. 
\end{conjecture}

Theorems \ref{thm:packingintro} and \ref{thm:packingintro2} are the first localized packing lemmas for subsets of groups of bounded $\VC_2$-dimension.  In the setting of hypergraphs,  a qualitative local packing lemma was proved by Chernikov and Towsner in \cite[Corollary 5.2]{Chernikov.2020}, and a distinct local packing result was obtained by Coregliano and Malliaris in \cite[Theorem 2.1]{Coregliano.2025}. The results in this paper complement strongly quantitative analogues in the hypergraph setting, to appear in forthcoming work by the first author.

 The remainder of the paper will be structured as follows.  Section \ref{ss:prelim} covers notation and preliminaries related to set systems, $\VC$-dimension, and quadratic factors.  We also give a proof of the failure of $\VC$-packing for quadratic atoms. In Section \ref{sec:arsep} we give a set of sufficient conditions for a set $A\subseteq \F_p^n$ to have $\VC_2(A)\geq k$.     In Section \ref{sec:packing}, we prove Theorem \ref{thm:packingintro} with the help of a linear cylinder regularity lemma, which will allow us to satisfy a uniformity condition that forms part of the set of sufficient conditions proved in Section \ref{sec:arsep}.  In Section \ref{sec:quasi}, we first prove Proposition \ref{exa:quadat} as a warm-up, after which we prove Theorem \ref{thm:packingintro2} from the structure theorem appearing in \cite{Sheats.2025} and some reasonably standard facts from linear algebra (including Lemma \ref{lem:subspace}, which to our knowledge is new).  The appendix contains a self-contained proof of the relevant cylinder regularity lemma (essentially a variation on \cite{Fox.2022}), followed by a proof of Lemma \ref{lem:subspace}.

\textbf{Acknowledgement.} The first author is partially supported by NSF CAREER Award DMS-2115518 and a Sloan Research Fellowship. The second author is supported by an Open Fellowship from the UK Engineering and Physical Sciences Research Council (EP/Z53352X/1). Both authors are deeply indebted to Maryanthe Malliaris, whose insightful questions about the authors' previous work \cite{Terry.2019,Terry.2021} led to the research in this paper and the forthcoming work \cite{Terry.2027}.

A detailed history of the paper and AI disclosure statement appears just before the bibliography.

\section{Preliminaries}\label{ss:prelim}
We first cover set systems, $\VC$-dimension, and $\VC_2$-dimension in groups in Subsection \ref{ss:vc}.  We then introduce quadratic factors in Subsection \ref{ss:atoms}, and prove Proposition \ref{prop:atom}, the failure of ordinary VC-packing for a single quadratic atom.  Finally Subsection \ref{ss:bool} sets notation related to boolean combinations and symmetric differences of subsets in $\F_p^n$.

\subsection{$\VC$-dimension and $\VC_2$-dimension}\label{ss:vc}
We cover here only the tools needed for the proofs that follow. We refer the reader to  \cite{Matousek.2009} for a more detailed introduction to $\VC$-theory in general, and to \cite{Alon.2019k3, Terry.2025, Terry.2025quo} for more background on $\VC$-dimension and $\VC_2$-dimension for subsets of groups.

Recall that a \emph{set system} is a pair $(X,\calF)$ where $\calF\subseteq \calP(X)$ is a collection of subsets of $X$. We call $X$ the \emph{ground set}.  We say a set system is \emph{finite} when $X$ is finite.

\begin{definition}[Shattering and $\VC$-dimension for set systems] \label{def:vc}
Let $(X,\calF)$ be a finite set system.  
\begin{enumerate}
\item We say that $\calF$ \emph{shatters a set} $T\subseteq X$ if for every $S\subseteq T$, there exists $F\in \calF$ such that $T\cap F=S$. 
\item The \emph{$\VC$-dimension of $(X,\calF)$}, denoted $\VC(X,\calF)$, is the largest integer $k\geq 0$ such that $\calF$ shatters a set of size $k$ in $X$.    
\item The \emph{dual of $(X,\calF)$} is the set system $(\calF, \calF^{dual})$ where $\calF^{dual}$ consists of the sets of the form $\calS_x=\{S\in \calF: x\in S\}$ for $x\in X$. 
\item The \emph{dual $\VC$-dimension of $(X,\calF)$} is $\VC(\calF,\calF^{dual})$.
\end{enumerate}
\end{definition}

The following relationship between the $\VC$-dimension of a set system  and its dual is well known (see, for example, \cite[Lemma 5.10]{Matousek.2009}).

\begin{fact}[$\VC$-dimension of the dual]\label{fact:dual}
Let $k\geq 1$, and let $(X,\calF)$ be a finite set system such that $\VC(X,\calF)<k$. Then $\VC(\calF,\calF^{dual})<2^{k}$, and this is best possible.
\end{fact}

We say a  collection of subsets $\calF\subseteq \calP(X)$ is \emph{$\delta$-separated} if for all distinct $F$ and $F'$ in $\calF$, $|F\Delta F'|>\delta|X|$.  Haussler's celebrated packing lemma, Theorem \ref{thm:haussler} below, gives a strong bound on the size of any $\delta$-separated collection when $\calF$ has bounded $\VC$-dimension \cite{Haussler.1995}.

\begin{theorem}[Haussler's packing lemma]\label{thm:haussler}
Let $k\geq 0$, let $\delta \in (0,1]$, and let $(X,\calF)$ be a finite set system such that $\calF$ is $\delta$-separated and such that $\VC(X,\calF)\leq k$. Then $\calF$ has size $|\calF|\leq (30/\delta)^k$.
\end{theorem}

We now turn to the setting of groups, specifically the elementary abelian $p$-groups that form the focus of this paper. We begin by defining the $\VC$-dimension of a subset of $G=\F_p^n$.

\begin{definition}[$\VC$-dimension for subsets of groups]\label{def:vcgroup}
Suppose $A\subseteq G=\F_p^n$. The \emph{$\VC$-dimension of $A$}, denoted $\VC(A)$, is the $\VC$-dimension of the set system with ground set $G$, and sets $\{A+g: g\in G\}$. 
\end{definition}

Our next task is to define a higher arity analogue of $\VC$-dimension called $\VC_2$-dimension, first defined by Shelah in \cite{Shelah.2017}.  To make its relationship to $\VC$-dimension clearer, we first observe that given $A\subseteq \F_p^n$ and $k\geq 1$, $\VC(A)\geq k$ if and only if there exist $a_1,\ldots, a_k\in G$ and $\{b_S: S\subseteq [k]\}\subseteq G$ such that $a_i+b_S\in A$ if and only if $i\in S$.  $\VC_2$-dimension is defined in an analogous way, shattering Cartesian products of sets, rather than sets of singletons.   

\begin{definition}[$\VC_2$-dimension for subsets of groups]\label{def:vc2group}
Let $A\subseteq G=\F_p^n$.
\begin{enumerate}
\item Given an integer $k\geq 1$, we say that $A$ \emph{has $\VC_2$-dimension at least $k$} if there exist $a_1,\ldots, a_k,b_1,\ldots, b_k\in G$ and $\{c_S: S\subseteq [k]^2\}\subseteq G$ such that $a_i+b_j+c_S\in A$ if and only if $(i,j)\in S$.
\item  The \emph{$\VC_2$-dimension of $A$}, denoted $\VC_2(A)$, is 
$$
\max\left(\left\{0\right\}\cup \left\{k\geq 1: A\text{ has $\VC_2$-dimension at least $k$} \right\}\right).
$$
\end{enumerate}
\end{definition}

It is easy to check that for any $A\subseteq \F_p^n$, $\VC_2(A)\leq \VC(A)$ (see \cite[Lemmas~3.1 and~3.2]{Terry.2021}). On the other hand, when $p$ is odd, the set $A=\{x\in \F_p^n: x\cdot x=0\}$ is an example of a set with $\VC_2$-dimension $1$, but whose $\VC$-dimension tends to infinity with $n$ (see \cite[Lemmas~3.9 and~3.10]{Terry.2025}, which originally appeared in an early version of \cite{Terry.2021}). This set $A$ is a prototypical example of what is known as an \emph{atom of a quadratic factor}, a notion which will be defined in detail in the next subsection.  

\subsection{Quadratic factors}\label{ss:atoms}

We begin with some terminology around quadratic factors.  

 \begin{definition}[Quadratic factor]
A \emph{quadratic factor} on $\F_p^n$ consists of a pair $\calB=(\calL,\calQ)$, where $\calL\subseteq \F_p^n$ is a set of linearly independent vectors, and $\calQ$ is a set of symmetric $n\times n$ matrices over $\F_p$.
\begin{enumerate}
\item The \emph{complexity of $\calB$} is $(\ell,q)$ where $\ell=|\calL|$ and $q=|\calQ|$.
\item When $q>0$, the \emph{rank of $\calB$} is  the minimum rank of any non-trivial linear combination of the elements of $\calQ$.  When $q=0$, then the rank of $\calB$ is $n$ by convention.
\end{enumerate}
When $\ell=0$, we refer to $\calB$ as simply $\calQ$.
\end{definition}

We next define the set of atoms corresponding to a quadratic factor.

\begin{definition}[Atoms of a quadratic factor]
Let $\calB=(\calL,\calQ)$ be a quadratic factor on $\F_p^n$ of complexity $(\ell,q)$.
\begin{enumerate}
\item If $\ell>0$, fix an enumeration $\calL=\{r_1,\ldots, r_{\ell}\}$, and for each $\abar=(a_1,\dots,a_\ell)\in \F_p^\ell$, define 
$$
L(\abar) =\{x\in G: x\cdot r_1=a_1,\ldots, x\cdot r_{\ell}=a_{\ell}\}.
$$
If $\ell=0$, let $L(0)=G$, where $0$ is the unique element of $\F_p^0$.
\item If $q>0$, fix an enumeration $\calQ=\{M_1,\ldots, M_q\}$, and for each $\bbar=(b_1,\ldots, b_q)\in \F_p^q$, define
\begin{align*} 
Q(\bbar)&=\{x\in G: x^TM_1x=b_1,\ldots, x^TM_qx=b_q\}.
\end{align*}
If $q=0$, let $Q(0)=G$.
\item Define $\At(\calB)=\{L(\abar)\cap Q(\bbar): \abar\bbar\in \F_p^{\ell}\times \F_p^q\}$ and refer to sets of the form $L(\abar)\cap Q(\bbar)$ as \emph{atoms of $\calB$}.
\end{enumerate}
\end{definition}

At several points in the paper we will use the following fact, which can be deduced by looking at the proof of \cite[Lemma 4.2]{Green.2007}.  

\begin{lemma}\label{lem:sizeofatoms}
Let $p>2$ be a prime, and let $(\calL,\calQ)$ be a quadratic factor on $\F_p^n$ of complexity $(\ell,q)$ and rank at least $r$.  Given $\abar\in \F_p^{\ell}$ and $\bbar\in \F_p^q$,  
$$
\left| |L(\abar)\cap Q(\bbar)| - p^{n-\ell-q}\right|\leq p^{n-r/2}.
$$
\end{lemma}

Using Lemma \ref{lem:sizeofatoms},  we can easily show that the quadratic atom $\{x\in \F_p^n: x\cdot x=0\}$ admits no efficient packing in the sense of Corollary \ref{thm:AFZ}. 

\begin{proposition}[$\VC$-packing fails for quadratic atoms]\label{prop:atom}
Let $p>2$ be a prime and let $n\geq 8$ be an integer. Let $G=\F_p^n$ and let 
\[A=\{x\in G: x\cdot x=0\}.\] 
For any $0<\delta<p^{-2}/2$, there exist elements $z_1,\ldots, z_n\in G$ such that $\{A-z_i: i\in [n]\}$ is $\delta$-separated.  
\end{proposition}
\begin{proof}
Let $n\geq 8$, and let $\{z_1,\ldots, z_n\}$ be any set of $n$ linearly independent vectors in $\F_p^n$. We show that the collection $\{A-z_i: i\in [n]\}$ is $\delta$-separated.  To this end, fix any $1\leq i\neq j\leq n$.  We shall show  that $|(A-z_i)\Delta (A-z_j)|\geq \delta |G|$. 

Let $\calL=\{r_1,r_2\}$ for $r_1=z_i$ and $r_2=z_j$,  and set $\calQ=\{I\}$, where $I$ is the $n\times n$ identity matrix.  Let $a_1=z_i\cdot z_i$, let $a_2=z_j\cdot z_j$, and set  
$$
J=\{((c_1,c_2),d)\in \F_p^2\times \F_p: d+2c_1+a_1=1\text{ and }d+2c_2+a_2=0\}.
$$
Note that if $x\in L(c_1,c_2)\cap Q(d)$ for some $((c_1,c_2),d)\in J$, then $x\in (A-z_j)\setminus (A-z_i)$.  It is easy to see that the set $J$ has exactly $p$ elements, and thus, using Lemma \ref{lem:sizeofatoms}, we have
\begin{align*}
|(A-z_i)\Delta (A-z_j)|\geq \sum_{((c_1,c_2),d)\in J}|L(c_1,c_2)\cap Q(d)|\geq |J|\left(p^{n-3}- p^{n/2}\right)&\geq p^{n-2}- p^{1+n/2},
\end{align*}
which is greater than $\delta |G|$, using $|G|=p^n$ and our assumptions on $\delta$ and $n$.
\end{proof}

\subsection{Boolean combinations and symmetric differences}\label{ss:bool}

Here we set notation for boolean combinations and symmetric differences.   Given sets $X,Y$, let $X^Y$ denote the set of functions from $Y$ to $X$.  When $X=\{0,1\}$, we write $2^Y$ to mean $X^Y$. Given a set $\Omega$ and $A\subseteq \Omega$, write $A^1=A$ and $A^0=\Omega\setminus A$. Whenever we use this notation, the ambient set $\Omega$ will be clear from context. We will use this notation mainly in the context of boolean combinations of an indexed family of sets.  For example, given $A_1,\ldots, A_s$, any atom of the Boolean algebra generated by $A_1,\ldots, A_s$  can be written as $\bigcap_{i=1}^s A_i^{\eta(i)}$, for some $\eta\in 2^{[s]}$.

\begin{definition} 
Let $\Omega$ be a finite set, let  $C\subseteq \Omega$ be a subset, and let $\calB=\{B_1,\ldots, B_t\}$ be a finite collection of subsets of $\Omega$. 
\begin{enumerate}
\item Define $\Bool_C(\calB)= \{C\cap \bigcap_{i=1}^t B_i^{\eta(i)}: \eta\in 2^{[t]}\}$.
\item Given an integer $k\leq t$, define $\Bool_C^{k}(\calB)=\{C\cap \bigcap_{i\in S}B_{i}^{\eta(i)}: \eta\in 2^S, S\in {[t]\choose k}\}$.
\item Define $\Sym_C(\calB)=\{(B_i\Delta B_j)\cap C: 1\leq i\neq j\leq t\}$ to be the set of all symmetric differences of pairs of elements of $\calB$, restricted to $C$.
\end{enumerate}
When $C=\Omega$, we may drop the subscript $C$.
\end{definition}
For later use, observe that in the notation above, $|\Bool_C(\calB)|\leq 2^t$, $|\Sym_C(\calB)|\leq {t\choose 2}$, and $|\Bool_C^{k}(\calB)|\leq 2^k{t\choose k}$.

\section{Sufficient conditions for large $\VC_2$-dimension}\label{sec:arsep}

The goal of this section is to prove Proposition \ref{prop:separation} below, which gives a sufficient condition for a set $A$ to have $\VC_2$-dimension at least $k$. This sufficient condition will be phrased in terms of the existence of a sufficiently large collection $\calS$ of translates of $A$, with the property that on some fixed coset $C$ the following two conditions hold:
\begin{enumerate}
\item $\calS$ is sufficiently separated on $C$;
\item certain boolean combinations from $\calS$ are sufficiently uniform on $C$. 
\end{enumerate}

Proving this proposition will require some preparatory work. The first lemma we shall prove, Lemma \ref{lem:packing} below, is a robust version of Haussler's packing lemma. To put it in context, note that the weaker conclusion that every element of $\Bool(\calS')$ is non-empty corresponds to ordinary Haussler packing.

\begin{lemma}[Robust Haussler packing]\label{lem:packing}
Let $k\geq 1$ be an integer, let $\rho\in (0,1/2]$, and set
$$
t=\lfloor(60/\rho)^{2^k}\rfloor+1\text{ and } \delta=\frac{\rho}{4{t\choose k}}.
$$ 
Let $\Omega$ be a finite set, and let $\calS$ be a collection of subsets of $\Omega$ of size $|\calS|=t$ such that for every $S\in \Sym(\calS)$, $|S|\geq \rho |\Omega|$. Then there exists a subcollection $\calS'\subseteq \calS$ of size $|\calS'|=k$ such that every element of $\Bool(\calS')$ has size at least $\delta |\Omega|$.   
\end{lemma}
\begin{proof}
Fix an enumeration $\calS=\{S_1,\ldots, S_t\}$.  Suppose towards a contradiction that  for all $I\in {[t] \choose k}$, there exists $\eta_{I}\in 2^{I}$ such that $|\bigcap_{i\in I}S_i^{\eta_{I}(i)}|\leq \delta |\Omega|$.  Let $X\subseteq \Omega$ be defined by setting
$$
X=\Omega\setminus \Bigg(\bigcup_{I\in {[t]\choose k}} \Bigg(\bigcap_{i\in I}S_i^{\eta_{I}(i)}\Bigg)\Bigg).
$$
Observe that
\begin{align}\label{al:omega}
|\Omega\setminus X|=\Bigg| \Bigg(\bigcup_{I\in {[t]\choose k}} \Bigg(\bigcap_{i\in I}S_i^{\eta_{I}(i)}\Bigg)\Bigg)\Bigg|\leq  {t\choose k}\delta|\Omega|=\frac{\rho}{4}|\Omega|,
\end{align}
where the last equality holds by definition of $\delta$.  We now consider the set system $(X, \calF)$ where $\calF=\{S_i\cap X: 1\leq i\leq t\}$.  Note that for all $1\leq i\neq j\leq t$, 
\begin{align*}
|(S_i\Delta S_j)\cap X|\geq |(S_i\Delta S_j)\cap \Omega|-|\Omega\setminus X|\geq \rho|\Omega|-\frac{\rho}{4}|\Omega|= \frac{3}{4}\rho |\Omega|>\frac{\rho|X|}{2},
\end{align*}
where the second inequality holds by (\ref{al:omega}).  Thus $\calF$ is $\rho/2$-separated on $X$.  

We next consider the dual set system of $\calF$. Writing $\calS_x=\{S_i: 1\leq i\leq t, x\in S_i\}$ for $x\in X$, we have
$$
\calF'=(\{S_1,\ldots, S_t\}, \{\calS_x : x\in X\}).
$$
We claim that no $k$-element subset of $\{S_1,\ldots, S_t\}$ is shattered in $\calF'$. Indeed, suppose there was $I\subseteq [t]$ with $|I|=k$ such that $\{S_i: i\in I\}$ was shattered in $\calF'$. That is, suppose that for all $J\subseteq I$, there is $x\in X$ such that $\{S_i: i\in J\}=\{S_i:i\in I, x\in S_i\}$. This would imply that for all $\sigma_I\in 2^I$, there is $x\in X\cap \bigcap_{i\in I} S_i^{\sigma_I(i)}$. But by construction, for each $I\in {[t]\choose k}$, $X\cap  \Big(\bigcap_{i\in I}S_i^{\eta_{I}(i)}\Big)=\emptyset$.  It follows that $\VC(\calF')<k$. 

By Fact \ref{fact:dual}, $\VC(\calF)<2^{k}$. Combining this with the earlier observation that $\calF$ is $\rho/2$-separated, Theorem \ref{thm:haussler} implies that $t\leq (60/\rho)^{2^k}$, a contradiction.
\end{proof}

We next define a local notion of uniformity which aligns with \cite[Definition  3.3]{Terry.2025} but has its precursor in \cite[Section 2]{Green.2005}.

\begin{definition}[Local uniformity norm]\label{def:locu2}
Let $G=\F_p^n$, let $H\leqslant G$ be a subgroup, and let $g\in G$. Given a function $f:H+g\rightarrow [-1,1]$, define
\[\|f\|_{U^2(g, H)}^4=\E_{\begin{subarray}{l} x_0,x_1\in H\end{subarray}}\E_{\begin{subarray}{l}y_0,y_1\in H+g\end{subarray}}\prod_{(i,j)\in \{0,1\}^2} f(x_{i}+y_{j}).\]
\end{definition}

In our context, we will be applying Definition \ref{def:locu2} to the balanced indicator function of a set, relativized to a fixed coset.  

\begin{definition}[Local uniformity on a coset]\label{def:unif} Given $A\subseteq G=\F_p^n$, $H\leqslant G$, and $g\in G$, we say that $A$ is \emph{$\e$-uniform on $H+g$} if $\|1_A-\alpha_{H+g}\|_{U^2(g,H)}\leq \epsilon$, where $\alpha_{H+g}=|A\cap(H+g)|/|H|$ is the density of $A$ on $H+g$.  
 \end{definition}

We next state a preliminary sufficient condition (Lemma \ref{lem:extip2}) for a set $A$ to have $\VC_2$-dimension at least $k$.  This lemma is already implicit in \cite{Terry.2021}, and it is also used in \cite{Terry.2025quo}.    For the sake of completeness we provide a deduction from the relevant counting lemma in \cite{Terry.2025} (which in turn is proved by a standard Cauchy-Schwarz argument).

\begin{lemma}[Sufficient condition for large $\VC_2$-dimension, I]\label{lem:extip2} 
Let $k\geq 1$ be an integer, let $\delta\in (0,1]$, and set 
\[\tau=\frac{\delta^{k2^{k^2}}}{(k+1)2^{k^2}}.\]

Let $G=\F_p^n$, let $H\leqslant G$ be a subspace, let $g\in G$ and let $A\subseteq G$. Suppose that $\{B_T: T\subseteq [k]\}$ is a collection of subsets of $H+g$ such that for all $T\subseteq [k]$,
\begin{enumerate}[label=\normalfont(\roman*)] 
\item $B_T$ is $\tau$-uniform on $H+g$, and
\item $|B_T|\geq \delta|H+g|$.
\end{enumerate}
Suppose further that $x_1,\dots,x_k\in G$ are such that for all $T\subseteq [k]$ and $j\in [k]$,
$$
B_T\subseteq A-x_j\text{ when }j\in T,\text{ and }B_T\cap (A-x_j)=\emptyset\text{ when }j\notin T.
$$ 

Then $A$ has $\VC_2$-dimension at least $k$.
\end{lemma}

\begin{proof}
Suppose first that $k=1$. In this case, choose $u_1=0$, $w_\emptyset \in B_\emptyset$, and $w_{\{(1,1)\}}\in B_{\{1\}}$. The assumptions give the required configuration.

Suppose now that $k>1$, and consider a family of arbitrary 1-bounded functions $(f_{i,S})_{i\in [k], S\subseteq [k]^2}$ defined on $H+g$. Let $m=k^2$, and identify the subsets of $[m]$ with the subsets of $[k]^2$ via a bijection $\phi:[m]\rightarrow [k]^2$, and for any $i\in [k+1,m]$ and any $U\subseteq [m]$, define $f_{i,\phi(U)}$ to be the constant function $1$.

We apply \cite[Lemma~3.8]{Terry.2025} with $m$, a linear factor whose atoms are the cosets of $H$, and labels $(a_1,a_2)$ corresponding to the cosets $H$ and $H+g$. With the above notation, 
\[\E_{\substack{u_i\in H: i\in [k]\\w_S\in H+g: S\subseteq [k]^2}} \prod_{i\in [k],S\subseteq [k]^2}f_{i,S}(u_i+w_S)=\E_{\substack{u_i\in H: i\in [m]\\w_U\in H+g: U\subseteq [m]}} \prod_{i\in [m],U\subseteq [m]}f_{i,\phi(U)}(u_i+w_U)\]
and thus it follows from \cite[Lemma~3.8]{Terry.2025} that
\begin{equation}\label{eq:telav}
\Big|\E_{\substack{u_i\in H: i\in [k]\\w_S\in H+g: S\subseteq [k]^2}} \prod_{i\in [k],S\subseteq [k]^2}f_{i,S}(u_i+w_S)\Big| \leq  \min_{i\in [k],S\subseteq [k]^2} \|f_{i,S}\|_{U^2(g,H)}.
\end{equation}
For each $T\subseteq[k]$, put $\beta_T=|B_T|/|H+g|$.
By assumption,
\[\beta_T\geq\delta \qquad\text{and}\qquad \|1_{B_T}-\beta_T\|_{U^2(g,H)}\le\tau.\]
For each $i\in[k]$ and $S\subseteq[k]^2$, write $S_i=\{j\in[k]:(i,j)\in S\}\subseteq [k]$.

We now outline a standard telescoping argument (see, for instance, \cite[Proof of Theorem 1.14]{Fox.2022}).  Starting with $f_{i,S}=1_{B_{S_i}}$ for all $i\in [k]$, $S\subseteq [k]^2$,  replace the product of $k2^{k^2}$ indicator functions on the left-hand side of \eqref{eq:telav} one at a time by the densities $\beta_{S_i}$. At each step, the difference between consecutive averages is an average over products of functions, one of which is of the form
$1_{B_{S_i}}-\beta_{S_i}$. It follows from \eqref{eq:telav} and the uniformity assumption that the error incurred with each replacement is at most $\tau$, so after $k2^{k^2}$ steps the total error is at most $k2^{k^2}\tau$.

We have thus shown that 
\[\E_{\substack{u_i\in H: i\in [k]\\w_S\in H+g: S\subseteq [k]^2}} \prod_{i\in [k],S\subseteq [k]^2}1_{B_{S_i}}(u_i+w_S) \geq  \prod_{i\in [k],S\subseteq [k]^2}\beta_{S_i} - k2^{k^2}\tau\geq \delta^{k2^{k^2}}-k2^{k^2}\tau =\frac{\delta^{k2^{k^2}}}{k+1}>0,\]
where the equality uses the definition of $\tau$.

There are thus $u_1,\dots,u_k\in H$ and for all $S\subseteq [k]^2$, $w_S\in H+g$ such that 
\begin{align}\label{al:uwb}
\text{for all $i\in [k]$ and $S\subseteq [k]^2$, $u_i+w_S\in B_{S_i}$. }
\end{align}
We now show that $\VC_2(A)\geq k$ by proving that for all $S\subseteq [k]^2$ and $i,j\in [k]$, $u_i+x_j+w_S\in A$ if and only if $(i,j)\in S$, equivalently, $u_i+w_S\in A-x_j$ if and only if $(i,j)\in S$.

Fix $S\subseteq [k]^2$ and $i,j\in [k]$.  Suppose first $(i,j)\in S$. By definition of $S_i$, we have $j\in S_i$, so $B_{S_i}\subseteq A-x_j$ by assumption. Consequently, by (\ref{al:uwb}), we have $u_i+w_S\in B_{S_i}\subseteq A-x_j$, as desired.  Suppose now $(i,j)\notin S$. Then $j\notin S_i$, so $B_{S_i}\cap (A-x_j)=\emptyset$. By (\ref{al:uwb}), we have $u_i+w_S\in B_{S_i}$, so this implies $u_i+w_S\notin A-x_j$, which finishes the proof.   
\end{proof}

We are now in a position to prove our main sufficient condition for a set to have large $\VC_2$-dimension.  

\begin{proposition}[Sufficient condition for large $\VC_2$-dimension, II]\label{prop:separation}
Let $k\geq 1$ be an integer, let $\rho\in (0,1/2]$, and let
$$
t=\lfloor (60/\rho)^{2^k}\rfloor+1\text{ , }\delta=\frac{\rho}{4{t\choose k}}\text{ , and }\tau=\frac{\delta^{k2^{k^2}}}{(k+1)2^{k^2}}.
$$
Let $G=\F_p^n$, let $H\leqslant G$, and let $A\subseteq G$. Suppose that there exist $g\in G$, $x_1,\ldots, x_t\in G$ such that for $\calS=\{(A-x_i)\cap (H+g): i\in [t]\}$ we have that
\begin{enumerate}[label=\normalfont(\roman*)] 
\item (uniformity) each $B\in \Bool_{H+g}^{ k}(\calS)$ is $\tau$-uniform on $H+g$;
\item (separation) for each $S\in \Sym_{H+g}(\calS)$, $|S\cap (H+g)|\geq \rho |H|$. 
\end{enumerate}

Then $A$ has $\VC_2$-dimension at least $k$.  
\end{proposition}
\begin{proof}
By Lemma \ref{lem:packing} applied to $\calS$ and $\Omega=H+g$, there exists a subcollection $\calS'$ of size $k$, say $\calS' =\{(A-x_{i_1})\cap (H+g), \ldots, (A-x_{i_k})\cap (H+g)\}$, such that every element of $\Bool_{H+g}(\calS')$ has size at least $\delta |H+g|$.  By assumption (i), every element of $\Bool_{H+g}^{k}(\calS')$ is also $\tau$-uniform on $H+g$.  Thus by Lemma \ref{lem:extip2}, $A$ has $\VC_2$-dimension at least $k$. 
\end{proof}

\section{Localized packing with polynomially many translates}\label{sec:packing}
In this section we prove Theorem \ref{thm:packingintro}. The proof will rely on the sufficient condition for large $\VC_2$-dimension from the preceding section (Proposition \ref{prop:separation}), along with a multicolor version of the $U^2$ cylinder regularity lemma. As is well understood at this point (see \cite{Fox.2022,Sheats.2025}), so-called ``cylinder partitions'', in which we do not require the partition to consist of cosets of a single subspace, but allow ourselves to refine by cosets of different subspaces, are quantitatively vastly superior. This will be essential to us in ultimately obtaining the bounds in Theorem \ref{thm:packingintro}. In our statements of these results, we will use  the following notion of a coset partition.  

\begin{definition}[Coset partitions]\label{def:cosetpart}
Let $G=\F_p^n$.
\begin{enumerate}
\item For a coset $C$ of a subgroup $H \leqslant G$, $\codim(C)$ denotes the codimension of $H$ in $G$.
\item A  \emph{coset partition of $G$} is a partition $\calC$ of $G$ such that for all $C\in \calC$, there exists a subspace $H_C\leqslant G$ such that $C$ is a coset of $H_C$.    
\item Given a coset partition $\calC$ of $G$, the \emph{codimension of $\calC$ in $G$} is defined to be 
$$
\codim(\calC)=\max\{\codim(C): C\in \calC\}.
$$  
\end{enumerate}
\end{definition}

We next state the cylinder regularity lemma we will use, Theorem \ref{thm:cyl} below. This result  is easily deduced from  \cite[Theorem 5.9]{Fox.2022}, although it is not identical to it. We include a self-contained proof in Appendix \ref{app:cylproof} for the sake of completeness, and to minimise the need for the reader to translate between different notations.

\begin{theorem}[$U^2$ cylinder regularity]\label{thm:cyl}
Let $p>2$ be a prime, let $s\geq 1$ be an integer and let $0<\eta  \leq 1$. Let $G=\F_p^n$ and suppose that $\calB=\{B_1,\ldots, B_s\}$ is a collection of subsets of $G$.  Then there exists a coset partition $\calC$ of $G$ and a set $\Omega\subseteq \calC$ such that
\begin{enumerate}[label=\normalfont(\roman*)] 
\item  $\codim(\calC)\leq \lceil s/ \eta^5\rceil$;
\item $\sum_{C\in \Omega}|C|\geq (1-\eta)|G|$;
\item for all $C\in \Omega$, every element of $\calB$ is $\eta$-uniform on $C$.
\end{enumerate}
\end{theorem}

We are now ready to prove Theorem \ref{thm:packingintro}.   

\vspace{2mm}

\begin{proofof}{Theorem \ref{thm:packingintro}}  Fix $k\geq 1$ and $0<\e<1/2$. Let  $t=t(k,\e)$ and $\tau=\tau(k,\e)$ be from Proposition \ref{prop:separation} applied with $\rho=\epsilon$, and set $M=t+1$.  Inspecting the form of $t$ and $\tau$ from Proposition \ref{prop:separation},  we find that there exist constants $a_k,b_k>1$, depending only on $k$, such that 
\begin{align}\label{al:tbound}
M\leq \e^{-a_k} \text{ and }\e^{b_k}\leq \tau<\e.
\end{align}
Define 
$$
L=\lceil (\e-\tau)^{-1}M\rceil +1 \text{ and }m=\lceil \tau^{-5}(2^kL^k)\rceil.
$$
In light of (\ref{al:tbound}), we see that there exist constants $r_k$ and $d_k$, depending only on $k$, such that $L\leq \e^{-r_k}$ and $m\leq \e^{-d_k}$. 

Let $G=\F_p^n$ and suppose that $A\subseteq G$ satisfies $\VC_2(A)< k$.  We now make a definition that will turn out to be crucial for the proof. Given a set $\calG\subseteq G$ and $x\in G$, we say that $x$ \emph{is $(\e,m)$-approximated by $\calG$} if there exists a coset partition $\calC$ of $G$ and a set $\Omega\subseteq \calC$ such that
\begin{itemize}
\item $\codim(\calC)\leq m$,
\item $|\bigcup_{C\in \Omega}C|\geq (1-\e)|G|$, and
\item for all $C\in \Omega$, there is $g_{x,C}\in \calG$ such that $|((A-x)\Delta (A-g_{x,C}))\cap C|\leq \e |C|$.   
\end{itemize}
Choose a sequence $(g_1,\ldots, g_R)\in G^R$ of maximal length such that for all $2\leq i\leq R$, $g_i$ is not $(\e,m)$-approximated by $\{g_1,\ldots, g_{i-1}\}$.  A maximal such sequence exists because the group $G$ is finite. By maximality of $(g_1,\ldots, g_R)$, we know that for every $x\in G$, there is a coset partition $\mathcal C_x$ of $G$ with $\codim(\calC_x)\leq m$, and a subcollection $\Omega_x\subseteq\mathcal C_x$, such that
\[\Big|\bigcup_{C\in\Omega_x}C\Big|\geq (1-\epsilon)|G|,\]
and such that for every $C\in\Omega_x$, there is some $j_{x,C}\in[R]$ satisfying
\[\left|\bigl((A-x)\triangle(A-g_{j_{x,C}})\bigr)\cap C \right|\leq \epsilon |C|.\]
To finish the proof, it suffices to show that $R< L$.   

Suppose towards a contradiction that $R\geq L$.  Let $\calB=\Bool^{k}(A-g_1,\ldots, A-g_L)$, and note that $|\calB|\leq {L\choose k}2^k \leq 2^kL^k$.   Apply Theorem \ref{thm:cyl} to $\calB$ with parameter $\tau$.  This yields a coset partition $\calC$ with $\codim(\calC)\leq m$, and a set $\Omega\subseteq \calC$ such that $|\bigcup_{C\in \Omega}C|\geq (1-\tau)|G|$, and such that every set in $\calB$ is $\tau$-uniform on every element of $\Omega$.  For each $2\leq i\leq L$, define
$$
\Sigma_i=\{C\in \calC: \text{ for all }1\leq j<i, |((A-g_i)\Delta(A-g_j))\cap C|>\e |C|\},
$$
i.e. $\Sigma_i$ is the set of those cosets in $\calC$ on which we are not able to approximate $A-g_i$ by any of the preceding translates. Now by construction, $g_i$ is not $(\e,m)$-approximated by $\{g_1,\ldots, g_{i-1}\}$, and thus for all $\Omega'\subseteq \calC$ with $|\bigcup_{C\in \Omega'}C|\geq (1-\epsilon)|G|$, $ \Omega'\cap \Sigma_i\neq \emptyset$. It follows that $|\bigcup_{C\in \Sigma_i}C|> \e|G|$, which in turn immediately implies that 
\begin{align}\label{al:in1}
\text{ for each $2\leq i\leq L$, }\left|\bigcup_{C\in \Omega\cap \Sigma_i}C\right|\geq \left|\bigcup_{C\in \Sigma_i}C\right|-\left|G\setminus \left(\bigcup_{C\in \Omega}C\right)\right|> (\e-\tau)|G|.
\end{align}
On the other hand, given $C\in \Omega$, define $N(C)=\{i\in [L]: i>1, C\in \Sigma_i\}$, i.e. $N(C)$ is the set of indices $i$ such that $A-g_i$ is not approximated on $C$ by an earlier translate. We claim that on any fixed $C\in \Omega$, $|N(C)|\leq M$. Fix $C\in \Omega$ and suppose towards a contradiction that $|N(C)|>M$. Then since we chose $M>t$, we can apply Proposition \ref{prop:separation} with $H+g=C$, $\rho=\e$, and an arbitrary selection of $t$ translates $A-g_{i_1},\dots, A-g_{i_t}$ with $i_1,\dots, i_t\in N(C)$. Indeed, recalling that every element of $\calB$ is $\tau$-uniform on $C$, and observing that by definition of $N(C)$, 
\[|((A-g_{i_a})\Delta(A-g_{i_b}))\cap C|> \epsilon |C|\]
for all $1\leq a\neq b\leq t$, Proposition \ref{prop:separation} implies that $\VC_2(A)\geq k$, contradicting our assumption on $A$.

Combining the above observations with (\ref{al:in1}), we find that 
\begin{align*}
(L-1)(\e-\tau)|G|< \sum_{i=2}^L\sum_{C\in \Omega \cap \Sigma_i}|C|\leq \sum_{C\in \Omega}|N(C)||C|\leq \sum_{C\in \Omega}M|C|\leq M|G|.
\end{align*}
Cancelling $|G|$ yields $(L-1)(\e-\tau)< M$. After rearranging, we thus have 
$$
L< M(\e-\tau)^{-1}+1,
$$
contradicting the definition of $L$. 
\end{proofof}

\section{Localized packing with subgroups of quasipolynomial index}\label{sec:quasi}

In this section we will prove Theorem \ref{thm:packingintro2} using the structure theorem appearing in \cite{Sheats.2025}.  As a warm-up, we first prove the localized packing lemma for the atom $\{x\in \F_p^n: x\cdot x=0\}$ stated as Proposition \ref{exa:quadat} in the introduction. The notation we use is deliberately more cumbersome than it strictly needs to be, in preparation for the more involved proof of Theorem \ref{thm:packingintro2}.
 
 \vspace{2mm}

\begin{proofof}{Proposition \ref{exa:quadat}}
Let $n\geq 3$, let $G=\F_p^n$, and let $A=\{x\in G: x\cdot x=0\}$. Let $v\in G\setminus\{0\}$ be such that $v\cdot v=0$ (its existence can be deduced from Lemma \ref{lem:sizeofatoms} and the fact that the identity matrix has rank $n$).  Fix any $h_0\in G$ such that $v\cdot h_0=1$, and set $h=h_0-\frac{1}{2}(h_0\cdot h_0)v$.  Note that $v\cdot h=1$ and $h\cdot h=0$.  Now set 
$$
\calG=\{tv: t\in \F_p\}\cup \{h+tv: t\in \F_p\}\subseteq G.
$$
Clearly $|\calG|\leq 2p$, and by construction $\calG\subseteq \langle v,h\rangle$.  We show $\calG$ has the desired properties.  Given $c\in G$, define $\Phi_c:G\rightarrow \F_p$ by setting
$$
\Phi_c(y)=y\cdot y+2c\cdot y.
$$

\begin{claim}\label{cl:0}
For all $a\in \F_p$ and $c\in G$, $\Phi_c^{-1}(a)\cap \calG\neq \emptyset$.
\end{claim}
\begin{proof}
Fix $c\in G$.  Let $\Psi_1,\Psi_2:\F_p\rightarrow \F_p$ be defined by setting $\Psi_1(t)=\Phi_c(tv)$ and $\Psi_2(t)=\Phi_c(h+tv)$.  Using that $v\cdot v=h\cdot h=0$ and $v\cdot h=1$, we have that for any $t\in \F_p$, $\Psi_1(t)=(2c\cdot v)t$ and $\Psi_2(t)=2(1+c\cdot v)t+2c\cdot h$.  Thus both $\Psi_1$ and $\Psi_2$ are (affine) linear maps with distinct slopes, and consequently at least one of them is injective. This implies that for all $a\in \F_p$, $\Phi_c^{-1}(a)\cap \calG\neq \emptyset$.
\end{proof}

Fix $x\in G$ and set $H_x=\langle v,h,x\rangle^{\perp}$. Clearly $H_x$ has codimension at most $3$.  Fix a coset $c+H_x$, and let $g$ be any element in $\Phi_c^{-1}(x\cdot x+2c\cdot x)\cap \calG$.  Our goal is to show that $(A-x)\cap (H_x+c)=(A-g)\cap (H_x+c)$, or equivalently, that for all $u\in H_x$, $u+c+x\in A$ if and only if $u+c+g\in A$. Since $A=\{x\in G: x\cdot x=0\}$, it suffices to show that for all $u\in H_x$, 
\begin{align}\label{al:mini}
(u+c+x)\cdot (u+c+x)=(u+c+g)\cdot (u+c+g).
\end{align}
Fix $u\in H_x$. Using that $u\in H_x\subseteq \langle x\rangle^{\perp}$, we have
$$
(u+c+x)\cdot (u+c+x)=u\cdot u+2u\cdot c+c\cdot c+(x\cdot x+2c\cdot x)=u\cdot u+2u\cdot c+c\cdot c+\Phi_c(x).
$$
On the other hand, using that $u\in H_x\subseteq \langle v,h\rangle^{\perp}$ and $g\in \calG\subseteq \langle v,h\rangle$, we have 
$$
(u+c+g)\cdot (u+c+g)=u\cdot u+2u\cdot c+c\cdot c+(g\cdot g+2c\cdot g)=u\cdot u+2u\cdot c+c\cdot c+\Phi_c(g).
$$
By our choice of $g$, $\Phi_c(g)=\Phi_c(x)$, so the preceding two displayed equations imply that (\ref{al:mini}) holds for this $u$. Since $u$ was an arbitrary element of $H_x$, this completes the proof.  
\end{proofof}

 \vspace{2mm}
 
Our next goal is to prove Theorem \ref{thm:packingintro2}.  We will need two lemmas about the rank of bilinear forms after restriction to a subspace.  To be precise, for a symmetric matrix $M$ and a subspace $V\leqslant G$, we shall write $M|_V$ for the restriction to $V\times V$ of the bilinear form $(x,y)\mapsto x^TMy$.

\begin{lemma}[Bounded rank drop on arbitrary bounded-codimension subspace]\label{lem:down}
Let $p>2$ be prime and let $r \geq 1$. Suppose that $\{M_1,\ldots, M_q\}$ is a quadratic factor on $G=\F_p^n$ of rank $r$, and $V\leqslant G$ is a subspace of  codimension $\ell$.  Then $\{M_1|_V,\ldots, M_q|_V\}$, considered as a quadratic factor on $V$, has rank at least $r-2\ell$. 
\end{lemma}

Lemma \ref{lem:down}, which is an easy exercise in linear algebra, follows immediately from \cite[Lemma 5.6]{Gowers.2012}, for example. The next lemma is less trivial, and allows us to find a small subspace preserving the rank of a set of matrices under restriction.

\begin{lemma}[No rank drop on special subspace]\label{lem:subspace}
Let $p>2$ be prime and let $r \geq 1$. Suppose that $\{M_1,\ldots, M_q\}$ is a quadratic factor on $G=\F_p^n$ of rank at least $r$. Then there is a subspace $W\leqslant G$ with
$\dim W\leq q+r$
such that
$\{M_1|_W,\ldots,M_q|_W\}$,
considered as a quadratic factor on $W$, has rank at least $r$.
 \end{lemma}
 
 We shall prove Lemma \ref{lem:subspace} in Appendix \ref{app:subspace}.    Let us now turn to the proof of Theorem \ref{thm:packingintro2}.
 
 \vspace{2mm}
 
\begin{proofof}{Theorem \ref{thm:packingintro2}}
Fix $k\geq 1$, $0<\e<1$, and a sufficiently large integer $n$.  Let $\rho(x)=4x+1$. Suppose that $A\subseteq \F_p^n$ has $\VC_2$-dimension less than $k$.  By \cite[Theorem 1.3]{Sheats.2025}, there exist integers 
$$
q\leq \log_p((\e/16)^{-2k-o(1)})\text{ and }\ell\leq \exp_p(\exp_p(\exp_p(O_{k,p,\rho}(\e^{-O(1)})))),
$$
a  quadratic factor $\calB=(\calL,\calQ)$ of complexity $(\ell,q)$ and rank at least $4\ell+4q+1$, and a set $\Omega\subseteq \At(\calB)$ such that $|A\Delta A^*|\leq \e^2|G|/16$, where $A^*=\bigcup_{B\in \Omega}B$.  To ease notation, for the duration of this proof, we let $V=L(0)$.  
 
We first deal with the case where $q=0$. In this case $A^*$ is a union of cosets of the subgroup $L(0)=V$. Set $r=p^{\ell}$ and choose $g_1,\ldots, g_r$ consisting of one representative of each coset of $V$.   For each $x\in G$, let $H_x=G$, and let $g_{x}$ be the element of $\{g_1,\ldots, g_r\}$ such that $g_{x}-x\in V$. Set $\Omega_x=\{G\}$.   Because $A^*$ is a union of cosets of $V$, $A^*-x=A^*-g_{x}$, and consequently 
\begin{align*}
|(A-x)\Delta (A-g_x)|&\leq |(A-x)\Delta (A^*-x)|+|(A^*-x)\Delta (A-g_x)|\\
&=|(A-x)\Delta (A^*-x)|+|(A^*-g_x)\Delta (A-g_x)|\\
&\leq \frac{\e^2|G|}{8}.
\end{align*}
This gives us the desired conclusion for the case $q=0$. 
 
Assume for the rest of the proof that $q>0$.  Since $V$ has codimension $\ell$ and $\calB$ has rank at least $4\ell+4q+1$, Lemma \ref{lem:down} implies $\calQ'=\{M_1|_{V},\ldots, M_q|_{V}\}$ has rank at least $4q+1$, considered as a quadratic factor on $V$.  Apply Lemma \ref{lem:subspace} with $r=4q+1$ to $\calB'=(\emptyset, \calQ')$ on $V$ to obtain a subspace $W\leqslant V$ of dimension at most $5q+1$ such that $\calQ''=\{M_1|_{W},\ldots, M_q|_{W}\}$, considered as a quadratic factor on $W$, has rank at least $4q+1$.  

For each $i\in [q]$ and $u,v\in G$, write $\beta_i(u,v)$ to denote the quantity $u^TM_iv$, and write $\beta(u,v)$ to denote the tuple $(\beta_1(u,v),\ldots, \beta_q(u,v))\in \F_p^q$.  For each $c\in G$, let $\Phi_c:G\rightarrow \F_p^q$ denote  the map 
$$
\Phi_c(x)=\beta(x,x)+2\beta(x,c).   
$$
\begin{claim}\label{cl:1}
For any $c\in G$ and $\abar\in \F_p^q$, 
\[\Big| |W\cap \Phi_c^{-1}(\abar)|-p^{\dim W-q}\Big|\leq p^{-1/2} p^{\dim W-q}.\] 
\end{claim}
\begin{proof}
Fix $c\in G$.  Fix a basis for $W$, and let $K_c=\langle M_1c,\ldots, M_qc\rangle^{\perp}\cap W$.   Let $\calL_c$ be a set of linearly independent vectors in $W$ such that $K_c=L_c(0)$ (here we use the dot product relative to the fixed basis for $W$).   Let $d_c=|\calL_c|$, and consider $(\calL_c,\calQ'')$ as a quadratic factor on $W$. By construction, it has complexity $(d_c,q)$ with $d_c\leq q$, and rank at least $4q+1$. Note that $W\cap \Phi_c^{-1}(\abar)$ is a union of $p^{d_c}$ atoms of $(\calL_c,\calQ'')$, and thus Lemma \ref{lem:sizeofatoms} implies 
\[\Big||W\cap \Phi_c^{-1}(\abar)|- p^{\dim W-q}\Big|\leq p^{d_c}p^{\dim W- 2q-1/2}\leq p^{-1/2} p^{\dim W-q},\]
as claimed.
\end{proof}
Recall we set $V=L(0)$, and let $U\leqslant G$ be any subspace satisfying $U\oplus V=G$, so that $\dim U=\codim V=\ell$. For all $x\in G$,  write $x=u_x+\ell_x$ where $u_x\in U$ and $\ell_x\in V=L(0)$,  and define  
$$
H_x=\{z\in G: \beta(z,\ell_x)=\beta(z,w)=0\text{ for all }w\in W\}.
$$
Then for each $C\in G/H_x$, choose any $v_{x,C}\in C$, and set 
$$
W_{x,C}=\{w\in W: \Phi_{v_{x,C}+u_x}(w)=\Phi_{v_{x,C}+u_x}(\ell_x) \}.
$$
(We shall see later in the proof that $W_{x,C}$ does not in fact depend on the choice of coset representative $v_{x,C}$.) By Claim \ref{cl:1}, for all $x\in G$ and $C\in G/H_x$,   we have 
\begin{align}\label{al:wsize}
\Big|| W_{x,C}|-p^{\dim W-q}\Big|\leq p^{-1/2} p^{\dim W-q}.
\end{align}
Note in particular that, since $\dim W$ is at least $4q+1$ (as $\calQ''$ has rank at least $4q+1$ in $W$), every set of the form $W_{x,C}$ is nonempty. 

Our next claim will imply that any translate of $A^*$ can be represented, locally, by a translate  by an element  from $U\oplus W$.  

\begin{claim}\label{cl:2}
For all $x\in G$, $C\in G/H_x$, and $w\in W_{x,C}$, $(A^*-x)\cap C=(A^*-u_x-w)\cap C$.
\end{claim}
\begin{proof}
Fix $x\in G$, a coset $C=H_x+v_{x,C}$, and any $w\in W_{x,C}$. Set $g=u_x+w$.  We show that $(A^*-x)\cap C=(A^*-g)\cap C$, or equivalently, that for all  $z\in C$, $z+x\in A^*$ if and only if $z+g\in A^*$. Since $A^*$ is a union of atoms of $\calB$, it suffices to show 
\begin{align}\label{al:claimgoal}
\text{ for all $z\in C$ and $(\abar,\bbar)\in \F_p^{\ell}\times \F_p^q$, $z+x\in L(\abar)\cap Q(\bbar)$ if and only if $z+g\in L(\abar)\cap Q(\bbar)$.}  
\end{align}
Fix $z\in C$. Since $z-v_{x,C}\in H_x$, we have $\beta(z,\ell_x)=\beta(v_{x,C},\ell_x)$ and $\beta(z,w)=\beta(v_{x,C},w)$.  Using $\beta(z,\ell_x)=\beta(v_{x,C},\ell_x)$, along with the fact that $x=u_x+\ell_x$, we have
\begin{align*}
\beta(z+x,z+x)&=\beta(z,z)+2\beta(u_x+\ell_x,z)+\beta(u_x+\ell_x,u_x+\ell_x)\\
&=\beta(z,z)+2\beta(u_x,z)+2\beta(\ell_x,v_{x,C})+\beta(u_x+\ell_x,u_x+\ell_x)\\
&=\text{$\beta(z,z)+2\beta(u_x,z)+\beta(u_x,u_x)+\Phi_{v_{x,C}+u_x}(\ell_x)$}.
\end{align*}
Similarly, using $\beta(z,w)=\beta(v_{x,C},w)$, along with the fact that $g=u_x+w$, we have
\begin{align*}
\beta(z+g,z+g)&=\beta(z,z)+2\beta(u_x+w,z)+\beta(u_x+w,u_x+w)\\
&=\beta(z,z)+2\beta(u_x,z)+2\beta(v_{x,C},w)+\beta(u_x+w,u_x+w)\\ 
&=\text{$\beta(z,z)+2\beta(u_x,z)+\beta(u_x,u_x)+\Phi_{v_{x,C}+u_x}(w)$}.
\end{align*}
By definition of $W_{x,C}$, and since $w\in W_{x,C}$, we have $\Phi_{v_{x,C}+u_x}(w)=\Phi_{v_{x,C}+u_x}(\ell_x)$. Therefore, the two observations above show that $\beta(z+x,z+x)=\beta(z+g,z+g)$.  On the other hand, given $\ell\in \calL$,  the definitions of $g$, $u_x$, and $\ell_x$ and the fact $w\in W\subseteq V= L(0)$, imply
\begin{align*}
(z+x)\cdot \ell=z\cdot \ell+u_x\cdot \ell =z\cdot \ell+u_x\cdot \ell+w\cdot \ell=(z+g)\cdot \ell.
\end{align*}
It follows that for any $(\abar,\bbar)\in \F_p^{\ell}\times \F_p^q$, $z+x\in L(\abar)\cap Q(\bbar)$ if and only if $z+g\in L(\abar)\cap Q(\bbar)$. Since $z$ was an arbitrary element of $C$, we conclude that (\ref{al:claimgoal}) holds, as desired.
 \end{proof}
 
 The remainder of the proof will show that since $A$ is approximated by $A^*$, and $A^*$ is locally approximated on cosets by translates from $U\oplus W$, so is $A$.

 For the rest of the proof, we fix $x\in G$.    Note that by definition of $H_x$ we have that for all $h\in H_x$ and $w\in W$,
 $$
 \Phi_{v_{x,C}+h+u_x}(w)= \beta(w,w)+2\beta(w,v_{x,C}+h+u_x)=\beta(w,w)+2\beta(w,v_{x,C}+u_x)=\Phi_{v_{x,C}+u_x}(w),
 $$
 and
  $$
 \Phi_{v_{x,C}+h+u_x}(\ell_x)= \beta(\ell_x,\ell_x)+2\beta(\ell_x,v_{x,C}+h+u_x)=\beta(\ell_x,\ell_x)+2\beta(\ell_x,v_{x,C}+u_x)=\Phi_{v_{x,C}+u_x}(\ell_x).
 $$
Consequently, the definition of $W_{x,C}$ is in fact independent of our choice of $v_{x,C}$ as a coset representative of $C$.  In other words, we have that for all $C\in G/H_x$ and $c\in C$, $W_{x,C}=\{w\in W: \Phi_{c+u_x}(w)=\Phi_{c+u_x}(\ell_x)\}$. Consequently, for any $y\in G$, using the notation $\calC_x(y)$ to denote the coset of $H_x$ containing $y$, we have
\begin{align}\label{al:wind}
\text{$W_{x,\calC_x(y)}=\{w\in W: \Phi_{y+u_x}(w)=\Phi_{y+u_x}(\ell_x)\}$.}
\end{align}

To ease notation, in what follows, we set $\Err=A\Delta A^*$, and recall that $|\Err|\leq \e^2|G|/16$.  For each $C\in G/H_x$, let  $\calG_{x,C}=\{u_x+w: w\in W_{x,C}\}$ and fix $g_{x,C}\in \calG_{x,C}$ such that 
$$
|(\Err-g_{x,C})\cap C|=\min_{g\in \calG_{x,C}}|(\Err-g)\cap C|.
$$ 
By definition of $\calG_{x,C}$, $g_{x,C}$ is of the form $u_x+w_{x,C}$ for some $w_{x,C}\in W_{x,C}$. Define
\begin{align*}
\Err_x^1&=\left\{C\in G/H_x: |((A^*-x)\Delta (A-x))\cap C|\geq \frac{\e|C|}{3}\right\}\\
\text{ and }\;\;\Err_x^2&=\left\{C\in G/H_x: |((A^*-g_{x,C})\Delta (A-g_{x,C}))\cap C|\geq \frac{\e|C|}{3}\right\}.
\end{align*}
Since $|\Err|=|\Err-x|=|(A-x)\Delta (A^*-x)|\leq \e^2|G|/16$, a standard averaging argument implies that $|\Err_x^1|\leq \e|G/H_x|/4$ for all $x\in G$. Controlling the size of $\Err_x^2$ will require more care. To start with, observe that
\begin{align*}
\sum_{C\in G/H_x} \sum_{g\in \calG_{x,C}}|(\Err-g)\cap C|&=\sum_{C\in G/H_x} \sum_{g\in \calG_{x,C}}|\Err\cap (C+g)|\\ 
&=\sum_{y\in \Err}|\{(C,g): C\in G/H_x, g\in \calG_{x,C}, y\in C+g\}|\\ 
&=\sum_{y\in \Err}|\{(C,g): C\in G/H_x, g\in \calG_{x,C}, C=\calC_x(y-g)\}|\\
&=\sum_{y\in \Err}|\{g\in G: g\in \calG_{x,\calC_x(y-g)}\}|.
\end{align*}
Given a fixed $y\in \Err$ and $g\in G$, we have $g\in \calG_{x,\calC_x(y-g)}$ if and only if $g=u_x+w$ for some $w\in W_{x,\calC_x(y-u_x-w)}$. Thus the above is equal to
\[\sum_{y\in \Err}|\{u_x+w:  w\in W_{x,\calC_x(y-u_x-w)}\}|=\sum_{y\in \Err}|\{w:  w\in W_{x,\calC_x(y-u_x-w)}\}|.\]
Now for a fixed $y\in \Err$, (\ref{al:wind}) implies that
\begin{align*}
\{w:  w\in W_{x,\calC_x(y-u_x-w)}\}&=\{w\in W: \Phi_{y-w}(w)=\Phi_{y-w}(\ell_x)\}\\
&=\{w\in W: -\beta(w,w)+2\beta(y+\ell_x,w)=\beta(\ell_x,\ell_x) +2\beta(\ell_x,y)\}\\
&=\{w\in W: \Phi_{-y-\ell_x}(w)=-\beta(\ell_x,\ell_x) -2\beta(\ell_x,y)\}.
\end{align*}

By Claim \ref{cl:1}, the right hand side has size at most $(1+p^{-1/2})p^{-q}|W|$. Combining the above with (\ref{al:wsize}), and recalling our choice of $g_{x,C}$, we have that 
\begin{align*}
\sum_{C\in G/H_x}|(\Err-g_{x,C})\cap C|&\leq \sum_{C\in G/H_x}\frac{1}{|W_{x,C}|} \sum_{g\in \calG_{x,C}}|(\Err-g)\cap C|\\
&\leq \frac{1}{p^{\dim W -q}(1-p^{-1/2})}\sum_{C\in G/H_x}  \sum_{g\in \calG_{x,C}}|(\Err-g)\cap C|\\
&= \frac{1}{p^{\dim W -q}(1-p^{-1/2})}\sum_{y\in \Err}|\{w:  w\in W_{x,\calC_x(y-u_x-w)}\}|\\
&\leq \frac{1+p^{-1/2}}{1-p^{-1/2}}|\Err|.
\end{align*}
Recalling that $|\Err|\leq \e^2|G|/16$, it follows that 
\[\sum_{C\in G/H_x}|(\Err-g_{x,C})\cap C|\leq \frac{\e^2(1+p^{-1/2})|G|}{16(1-p^{-1/2})}< \frac{\e^2|G|}{4}.\]
By a standard averaging argument, the above inequality and the definition of $\Err_x^2$ imply that $|\Err_x^2|<3\e|G/H_x|/4$. 

It remains to show that $A-x$ is well approximated by a shift from $\calG_{x,C}$ on any coset $C$ not included in $\Err_x^1\cup \Err_x^2$. To this end, define
$$
\Omega_x=(G/H_x)\setminus (\Err_x^1\cup \Err_x^2),
$$
and note that $|\Omega_x|>(1-\e)|G/H_x|$. Then for all $x\in G$ and $C\in \Omega_x$, 
\begin{align*}
|((A-x)\Delta (A-g_{x,C}))\cap C|\leq |((A-x)&\Delta (A^*-x))\cap C|+|((A^*-x)\Delta (A^*-g_{x,C}))\cap C|\\
&+|((A-g_{x,C})\Delta (A^*-g_{x,C}))\cap C|,
\end{align*}
which, by definition of $\Err_x^1$, $\Err_x^2$, and Claim \ref{cl:2}, is at most
\[\frac{\e|C|}{3}+0+\frac{\e|C|}{3}<\e |C|.\]
By construction each $H_x$ has codimension at most 
$$
D\leq q(\dim W+1)\leq 5q^2+2q\leq 8\log^2_p(\e^{-2k-o(1)}),
$$
and each $g_{x,C}$ comes from $U\oplus W$, a set of size at most 
$$ 
p^{\ell+5q+1}\leq \exp_p(\exp_p(\exp_p(\exp_p(O_{p,k,\rho}(\e^{-O(1)}))))).
$$
This completes the proof.
\end{proofof}

\appendix

\section{Proof of Theorem \ref{thm:cyl}}\label{app:cylproof}

In this appendix we prove Theorem \ref{thm:cyl}.  We first state more definitions related to Definition \ref{def:cosetpart}.

\begin{definition}[Coset partitions]\label{def:cosetpart2}
Let $G=\F_p^n$, and let $L$ be a coset of a subgroup $H \leqslant G$.
\begin{enumerate} 
\item A  \emph{coset partition of $L$} is a partition $\calC$ of $L$ such that for all $C\in \calC$, there exists a subspace $H_C\leqslant H$ such that $C$ is a coset of $H_C$.  
\item Given a coset partition $\calC$ of $L$, the \emph{codimension of $\calC$ in $L$} is defined to be 
$$
\codim(\calC;L)=\max\{\codim(H_C; H): C\in \calC\}.
$$ 
\item Given $x\in L$ and a coset partition $\calC$ of $L$, $\calC(x)$ denotes the element $C$ of $\calC$ containing $x$.
\item Given two coset partitions $\calC$ and $\calC'$ of $L$, we say $\calC'$ \emph{refines $\calC$}, denoted $\calC'\preceq \calC$, if $\calC'$ refines $\calC$ as a partition of $L$. 
\end{enumerate}
\end{definition}

The proof of our cylinder regularity lemma will gradually increment the $L^2$ energy relative to a coset partition, a quantity which we refer to as the \emph{index}.

\begin{definition}[Index of a set with respect to a coset partition]\label{def:indcos}
Let $G=\F_p^n$, let $L$ be a coset of a subgroup of $G$, and let $\calC$ be a coset partition of $L$. Define the \emph{index of a subset $A\subseteq G$ with respect to  $\calC$ in $L$} by 
\[
\ind_L(A;\calC)=\frac{1}{|L|}\sum_{x\in L}\Bigg(\frac{|A\cap \calC(x)|}{|\calC(x)|}\Bigg)^2=\frac{1}{|L|}\sum_{C\in \calC}|C|\Bigg(\frac{|A\cap C|}{|C|}\Bigg)^2.
\]
\end{definition}

We now extend these notations to collections of subsets of a group.
 
\begin{definition}[Index of a collection with respect to a coset partition]\label{def:indcos2}
Let $G=\F_p^n$,  let $L$ be a coset of a subgroup of $G$, and let $\calC$ be a coset partition of $L$. Given a collection $\calB=\{B_1,\ldots, B_s\}$ of subsets of $G$, define the \emph{index of $\calB$ with respect to $\calC$ in $L$} by
\[\ind_L(\calB;\calC)=\frac{1}{|L|}\sum_{x\in L}\sum_{i=1}^s \Bigg(\frac{|B_i\cap \calC(x)|}{|\calC(x)|}\Bigg)^2=\sum_{i=1}^s\ind_L(B_i;\calC).\]
\end{definition}

The following lemma is immediate from Proposition 5.2 parts (1) and (2) in \cite{Fox.2022}. 

\begin{lemma}[Standard properties of index]\label{lem:mon}
Let $G=\F_p^n$, let $L$ be a coset of a subgroup of $G$,  let $\calC$ be a coset partition of $L$, and let $\calB=\{B_1,\ldots, B_s\}$ be a collection of subsets of $G$. Then the following hold.
\begin{enumerate}
\item $0\leq \ind_L(\calB; \calC)\leq s$.
\item If  $\calC'$ is another coset partition of $L$ satisfying $\calC'\preceq \calC$, then
\begin{align*} 
\ind_L(\calB;\calC')-\ind_L(\calB; \calC)=\frac{1}{|L|}\sum_{i=1}^s\sum_{C\in \calC}\left(\sum_{C'\in \calC', C'\subseteq C}|C'|\Bigg(\frac{|B_i\cap C|}{|C|}-\frac{|B_i\cap C'|}{|C'|}\Bigg)^2\right)\geq 0. 
\end{align*}
\end{enumerate}
\end{lemma}

We will use the fact that large uniformity leads to a density increment on a coset of a codimension-1 subspace.  

\begin{lemma}[Lack of uniformity implies energy increment]\label{lem:densincr}
Let $G=\F_p^n$, let $H \leqslant G$ be a subgroup, and let $g\in G$. Let $\eta\in (0,1]$, and let $A\subseteq G$ be a set which is not $\eta$-uniform on $L=H+g$. Then there exists a coset partition $\calC$ of $L$ of codimension $1$ such that  
$$
\ind_L(A;\calC)-\ind_L(A; \{L\})\geq \eta^4.
$$
\end{lemma}

Lemma \ref{lem:densincr} consists of standard Fourier arguments, and is essentially the same as \cite[Proposition 5.2(3)]{Fox.2022}. We include a proof here as we have chosen a slightly different notion of uniformity in this paper.  We will need a few preliminaries. First, we let $\omega=e^{2\pi i/p}$, and for a function $f:G=\mathbb F_p^n\to\C$, we define its Fourier transform by
\[
\widehat f(t):=\E_{x\in G}
f(x)\omega^{-t\cdot x}.
\]
With this normalization, the Fourier inversion formula takes the form
$$
f(x)=\sum_{t\in\widehat{G}}
\widehat f(t)\omega^{t\cdot x},
$$
and Parseval's identity is
$$
\E_{x\in G}|f(x)|^2=\sum_{t\in\widehat{G}}
|\widehat f(t)|^2.
$$
We are now ready to prove Lemma \ref{lem:densincr}. 

\vspace{2mm}

\begin{proofof}{Lemma \ref{lem:densincr}} Set $\alpha_{H+g}=|A\cap(H+g)|/|H|$, and define
$f:H\to[-1,1]$ by $f(x)=1_A(x+g)-\alpha_{H+g}$. Note that $\E_{x\in H}f(x)=0$, and that 
\[\|1_A-\alpha_{H+g}\|_{U^2(g,H)}=\|f\|_{U^2(0,H)}= \|f\|_{U^2(H)},\]
where the rightmost term is simply the standard $U^2$ norm defined on $H$. It is well known that 
\[\|f\|_{U^2(H)}^4=\sum_{t\in \widehat{H}}|\widehat{f}(t)|^4 \leq \sup_{t\in \widehat{H}}|\widehat{f}(t)|^2 \sum_{t\in \widehat{H}}|\widehat{f}(t)|^2,\]
so by Parseval's identity the failure of $\eta$-uniformity gives $t\in \widehat{H}$ such that 
$|\widehat{f}(t)|\geq \eta^2$.

Let $H'=\langle t\rangle^{\perp}$ in $H$, and let $\calC$ be the partition of $L=H+g$ into cosets of $H'$.  Note we may enumerate these cosets as $H_0+g,\ldots, H_{p-1}+g$, where $H_u=\{x\in H: x\cdot t=u\}$ for each $u\in \{0,\ldots, p-1\}$. We then have
\[\eta^2\leq |\widehat{f}(t)|=|\E_{x\in H}f(x)\omega^{-t\cdot x}|=|\E_{u\in \F_p}\E_{x\in H_u}f(x)\omega^{-u}|\leq \E_{u\in \F_p}|\E_{x\in H_u}f(x)| .\]
But by definition of $f$, the expectation in $x$ is
\[\E_{x\in H_u}f(x)=\E_{x\in H_u}(1_A(x+g)-\alpha_{H+g})=\alpha_{H_u+g}-\alpha_{H+g},\]
and so
\[\eta^2\leq\E_{u\in \F_p} |\alpha_{H_u+g}-\alpha_{H+g}|.\]
On the other hand, by Lemma \ref{lem:mon},
\[\ind_L(A;\calC)-\ind_L(A;\{L\})=\frac{1}{|H|}\sum_{u=0}^{p-1}|H_u|(\alpha_{H_u+g}-\alpha_{H+g})^2=\E_{u\in \F_p}(\alpha_{H_u+g}-\alpha_{H+g})^2\]
which, by Cauchy-Schwarz, is at least
\[\left(\E_{u\in \F_p}|\alpha_{H_u+g}-\alpha_{H+g}|\right)^2\geq \eta^4,\]
as desired.
\end{proofof}
 
 \vspace{2mm}
 
 We are now ready to prove  Theorem \ref{thm:cyl}.
 
 \vspace{2mm}

\begin{proofof}{Theorem \ref{thm:cyl}}
We proceed via an energy increment argument.

\noindent{\bf Step $0$:} Let $\calC_0=\{G\}$.  If every element of $\calB$ is $\eta$-uniform with respect to $G$, set $\calC=\calC_0$ and $\Omega=\calC$, and end the construction. Otherwise go to the next step.

\noindent{\bf Step $i+1$:} Suppose $0\leq i<\lceil  s/ \eta^5 \rceil+1$ and assume by induction that we have defined a coset partition $\calC_i$ of $G$ satisfying
\begin{itemize}
\item $\ind_G(\calB; \calC_i)\geq i \eta^5$;
\item $\codim(\calC_i)\leq i$;
\item the set $\Sigma_i=\{C\in \calC_i: \text{some element $B_{j_C}\in \calB$ is not $\eta$-uniform on $C$}\}$ has the property that $\sum_{C\in \Sigma_i}|C| \geq \eta |G|$. 
\end{itemize}
For each $C=H_C+g_C\in \Sigma_i$, Lemma \ref{lem:densincr} implies there is a coset partition $\calD_C$ of $C$ satisfying $\codim(\calD_C; C)=1$ such that $\ind_C(B_{j_C}; \calD_C)- \ind_C(B_{j_C}; \{C\})\geq \eta^4$.  Combined with Definition \ref{def:indcos2}, this implies $\ind_C(\calB; \calD_C)- \ind_C(\calB; \{C\})\geq \eta^4$. 

Now let 
$$
\calC_{i+1}=(\calC_i\setminus  \Sigma_i)\cup \bigcup_{C\in \Sigma_i }\calD_C
$$
By construction and our induction hypothesis, $\codim(\calC_{i+1})\leq 1+\codim(\calC_i)\leq i+1$. We can express    $\ind_G(\calB; \calC_{i+1})-\ind_G(\calB; \calC_i)$ as
\begin{align*}
\Big(\sum_{C\in \Sigma_i}\sum_{C'\in \calC_{i+1}, C'\subseteq C}\frac{|C'|}{|G|}&\ind_{C'}(\calB; \{C'\})\Big)-\Big(\sum_{C\in \Sigma_i}\frac{|C|}{|G|}\ind_C(\calB; \{C\})\Big)\\
=&\sum_{C\in \Sigma_i}\frac{|C|}{|G|}\Big(\sum_{C'\in \calC_{i+1}, C'\subseteq C}\frac{|C'|}{|C|}\ind_{C'}(\calB; \{C'\}) - \ind_C(\calB; \{C\})\Big)\\
=&\sum_{C\in \Sigma_i}\frac{|C|}{|G|}\Big(\ind_C(\calB; \calD_C)   - \ind_C(\calB; \{C\})\Big),
\end{align*}
which by our assumption on $\Sigma_i$, and our construction of $\calC_{i+1}$,  is bounded below by
\[\eta^4\sum_{C\in \Sigma_i}\frac{|C|}{|G|}\geq \eta^5.\]
Consequently, $\ind_G(\calB; \calC_{i+1})\geq \ind_G(\calB; \calC_i)+\eta^5\geq (i+1)\eta^5$, where the second inequality uses the induction hypothesis.  

Now if $\calC_{i+1}$ satisfies the desired conclusions, then set $\calC=\calC_{i+1}$ and end the construction. Otherwise we can go on to the next step. Clearly this proceeds at most $\lceil s/\eta^5 \rceil$ times, which finishes the proof.
\end{proofof}

\section{Proof of Lemma \ref{lem:subspace}}\label{app:subspace}
 
We restate Lemma \ref{lem:subspace} for the convenience of the reader.

\begin{proposition}
Let $p>2$ be prime and let $r \geq 1$. Suppose that $\{M_1,\ldots, M_q\}$ is a quadratic factor on $V=\F_p^n$ of rank at least $r$. Then there is a subspace $H\leqslant V$ with
$\dim H\leq q+r$
such that
$\{M_1|_H,\ldots,M_q|_H\}$,
considered as a quadratic factor on $H$, has rank at least $r$.
\end{proposition}

\begin{proof}
If $n\leq q+r$, we may simply let $H=V$ and finish the proof.

Assume now that $n>q+r$.   Let $\calQ=\{M_1,\ldots,M_q\}$. The rank hypothesis implies that the $M_i$ are linearly independent, so $W=\mathrm{span}(\calQ)$ has $\dim W=q$.  We shall find a hyperplane $H\leqslant V$  such that
\[\operatorname{rank}(M|_H)\geq r\]
for every non-zero $M\in W$.

Let $V^*$ denote the space of linear maps from $V$ to $\F_p$. For $M\in W$, let $R(M)\leqslant V^*$ denote its row space. Consider
those non-zero $M\in W$ whose rank is either $r$ or $r+1$.
Scalar multiples of a fixed $M$ have the same row space, and the number of
one-dimensional subspaces of $ W$ is
\[\frac{p^q-1}{p-1}.\]
Moreover, each relevant row space contains at most $p^{r+1}$ elements.
Consequently,
\[\Big|\bigcup_{\substack{M\in W\setminus\{0\}\\ \mathrm{rank}(M)\leq r+1}}
R(M)\Big| \leq \frac{p^q-1}{p-1}\,p^{r+1} <p^{q+r+1} \leq p^n.\]
We may therefore choose a non-zero $\lambda\in V^*$ outside this union,
and set $H=\ker\lambda$. Since $\lambda$ is non-zero, $\dim(H)= n-1$. We verify that $H$ has the required property. 

First suppose that $\mathrm{rank}(M)\geq r+2$. By Lemma \ref{lem:down}, $\mathrm{rank}(M|_H)\geq \mathrm{rank}(M)-2\geq r$.

Now suppose that $\mathrm{rank}(M)$ is either $r$ or $r+1$.
We claim that restriction to $H$ does not lower the rank of $M$. To see this, consider
\[\mathrm{rad}(M|_H)=\{x\in H:x^{T}My=0\text{ for every }y\in H\}.\]
If $x\in\mathrm{rad}(M|_H)$, then the functional $x^{T}M$ vanishes on $H=\ker\lambda$. It is therefore a scalar
multiple of $\lambda$, i.e. $x^{T}M=c\lambda$ for some $c\in\F_p$. Since $x^{T}M\in R(M)$, whereas $\lambda\notin R(M)$, we must have $c=0$. Thus $x^{T}M=0$, and symmetry gives $Mx=0$. It follows that
\[\mathrm{rad}(M|_H)=H\cap\ker M.\]

Elementary linear algebra shows that the row space of $M$ is the annihilator of its kernel, i.e.
$R(M)=(\ker M)^\circ=\{S\in V^*: Sx=0  \mbox{ for all }x \in \ker M\}$. Since $\lambda\notin R(M)$, the restriction of $\lambda$ to $\ker M$ is non-zero. Therefore
\[\dim(H\cap\ker M)=\dim(\ker M)-1.\]
It follows that
\[\mathrm{rank}(M|_H)=\dim H-\dim(\mathrm{rad}(M|_H))=(n-1)-(\dim(\ker M)-1)=\mathrm{rank}(M),\]
and so every non-zero element of $ W$ still has rank at least $r$ after restriction to $H$. In particular, the restrictions of $M_1,\ldots,M_q$ remain linearly independent. We may therefore repeat the argument, starting with $V=H$, reducing the ambient dimension by one at each step, until we
obtain a subspace $H\leqslant \F_p^n$ with
\[\dim H\leq q+r.\]
Every non-zero member of $ W$ restricts to a form of rank at least $r$ on $H$, as required.
\end{proof}

\textbf{History of the paper and AI involvement.} The mathematical ideas in this paper date back to the summer of 2022. Over the subsequent months, the authors developed a proof of a local packing result with wowzer-type bounds, which were a consequence of an application of a strong linear regularity lemma. In early 2026, the authors simplified the argument using a linear cylinder decomposition, which reduced the bound to a tower. In June 2026, ChatGPT Pro 5.5 first observed that the proof ingredients could be rearranged to give a double exponential bound, and that a slightly weaker version of our separation lemma was in fact sufficient. In this and earlier iterations, the packing statement was qualitatively weaker, requiring different approximating translates on different parts of a fixed coset partition. Two days later, in June 2026, in parallel work by the first author on hypergraphs, ChatGPT Pro 5.5 suggested a further rearrangement of the proof ingredients which produced the qualitatively stronger result presented in this paper as Theorem \ref{thm:packingintro}, and which resulted in a single exponential bound.  

It was known to the authors since 2022 that a localized packing lemma holds for unions of quadratic atoms, and thus that one can obtain a version of the proof of Theorem \ref{thm:packingintro2}.  Indeed, these observations gave rise to the notion of localized packing appearing in this paper.  The eventual proofs of Proposition \ref{exa:quadat} and Theorem \ref{thm:packingintro2} appearing in this paper have benefited from input by ChatGPT Pro 5.6, which streamlined and optimised the quantitative aspects of the arguments.  

ChatGPT Pro 6 was also used for general proofreading. All input from ChatGPT has been thoroughly checked and reworked by the authors.


\providecommand{\bysame}{\leavevmode\hbox to3em{\hrulefill}\thinspace}
\providecommand{\MR}{\relax\ifhmode\unskip\space\fi MR }
\providecommand{\MRhref}[2]{%
  \href{http://www.ams.org/mathscinet-getitem?mr=#1}{#2}
}
\providecommand{\href}[2]{#2}

\end{document}